\documentclass[11pt,reno]{amsart}  
\usepackage{amssymb}
\usepackage{amsmath}
\usepackage{enumitem}
\usepackage{mathrsfs}
\usepackage[english]{babel}
\usepackage[all]{xy}
\usepackage{hyperref}
\usepackage{marginnote}
\usepackage{stmaryrd}
\usepackage{fullpage}

\RequirePackage{mathrsfs}

\newtheorem{theorem}{Theorem}

\theoremstyle{definition}
\newtheorem*{ack}{Acknowledgements}
\newtheorem{definition}{Definition}

\DeclareMathOperator{\Spec}{Spec}

\title{A mixed fibration theorem for Hilbert irreducibility on non-proper varieties}
\author{Cedric Luger}
\address{Cedric Luger\\
Institut f\"{u}r Mathematik\\
Johannes Gutenberg-Universit\"{a}t Mainz\\
Staudingerweg 9, 55099 Mainz\\
Germany.}
\email{celuger@uni-mainz.de}

\subjclass[2020]
{
    14G99  
    (14G05,  
    14G40)}  

\keywords{
	Integral points,
	ramified covers,
	Hilbert irreducibility}

\begin{document}

\begin{abstract}
    We prove that the weak Hilbert property ascends along a morphism
    of varieties over an arbitrary field of characteristic zero, under suitable assumptions. 
\end{abstract}

\maketitle

\thispagestyle{empty}

A variety $X$ over a field $k$ is said to have the \emph{Hilbert property over $k$} if its $k$-rational points are not thin (cf. \cite[\S 3]{SerreTopicsGalois}).
\begin{definition}\label{Def:HP}
    Let $k$ be a field of characteristic zero,
    let $X$ be an integral variety over $k$,
    and let $\Sigma\subseteq X(k)$ be a subset.
    We say that $\Sigma$ is \emph{thin} in $X$ if there exists an integer $n\geq 0$
    and a finite collection $(\pi_i \colon Y_i \to X)_{i=1}^n$ of finite surjective morphisms of integral varieties of degree at least~$2$ such that
    \[
        \Sigma \setminus \bigcup_{i=1}^n \pi_i(Y_i(k)) \subseteq X(k)
    \]
    is not Zariski-dense in $X$. We say that $X$ satisfies the \emph{Hilbert property (over~$k$)} if $X(k)$ is not thin in $X$.
\end{definition}
Note that, in Definition~\ref{Def:HP}, we may always assume the $Y_i$ to be normal (by replacing them with a normalization) and geometrically integral (as they would not have a $k$-rational point otherwise).
The Hilbert property is usually defined using generically finite covers;
this gives the same notion since we may replace a generically finite cover $Y_i \to X$ by the relative normalization of $X$ in $Y_i$.

Serre asked whether the product of two varieties with the Hilbert property satisfies the Hilbert property (see the problem stated in \cite[\S 3.1]{SerreTopicsGalois}).
This question was answered positively by Bary-Soroker, Fehm, and Petersen \cite[Corollary~3.4]{BSFP14}.
More generally, Bary-Soroker--Fehm--Petersen prove a fibration theorem \cite[Theorem~1.1]{BSFP14} for the Hilbert property: If $X \to S$ is a dominant morphism of $k$-varieties such that $S$ and all fibres $X_s$ with $s\in S(k)$ satisfy the Hilbert property, then $X$ satisfies the Hilbert property.

Corvaja--Zannier \cite[Theorem~1.4]{CZHP} showed that,
if $k$ is a number field and $X$ is a normal integral projective $k$-variety,
then the presence of a non-trivial finite étale cover of $X$ forces $X(k)$ to be thin.
In light of this, they introduce the \emph{weak Hilbert property} \cite[Definition~1.2]{CDJLZ},
in which only ramified covers are considered.

\begin{definition}
    Let $X$ be a normal integral variety over a field $k$ of characteristic zero.
    A \emph{ramified cover} of $X$ is a finite surjective morphism $Y \to X$ of normal integral varieties over $k$ that is not unramified.
\end{definition}

\begin{definition}\label{Def:WHP}
    Let $k$ be a field of characteristic zero,
    let $X$ be a normal integral variety over $k$,
    and let $\Sigma\subseteq X(k)$ be a subset.
    We say that $\Sigma$ is \emph{strongly thin in $X$} if there exists an integer $n\geq 0$
    and a finite collection $(\pi_i \colon Y_i \to X)_{i=1}^n$ of ramified covers of normal integral varieties such that
    \[
        \Sigma \setminus \bigcup_{i=1}^n \pi_i(Y_i(k)) \subseteq X(k)
    \]
    is not Zariski-dense in $X$.
    We say that $X$ satisfies the \emph{weak Hilbert property (over $k$)} if $X(k)$ is not strongly thin in $X$.
\end{definition}

Definition~\ref{Def:WHP} is restricted to normal varieties due to the fact that in this case a cover $Y \to X$ is either ramified or étale (see \cite[Lemma~2.3]{CDJLZ}). Moreover, a ramified cover is of degree at least~$2$
(as otherwise it would be an isomorphism by Zariski's main theorem),
which shows that a strongly thin set is thin.

It has been shown by Javanpeykar and Wittenberg (see \cite[Theorem~1.9]{CDJLZ}) that
the product of two smooth proper varieties over a finitely generated field $k$ of
characteristic zero with the weak Hilbert property over $k$ satisfies the weak Hilbert property over $k$.
(Without the properness assumption, a similar statement holds by considering near-integral points, cf. \cite[Theorem~1.4]{LugerProducts}).
The proof of the product theorem \cite[Theorem~1.9]{CDJLZ} does allow to consider subsets of rational points that are not necessarily given by a product of sets of rational points (see \cite[Theorem~4.4]{LugerProducts} for a precise statement).
However, it still requires to work with a product of varieties, instead of a more general fibration as in \cite[Theorem~1.1]{BSFP14}.

A fibration theorem similar to the result of Bary-Soroker--Fehm--Petersen for the weak Hilbert
property was established in \cite[Theorem~1.3]{JavanpeykarNefTangentBundle}, where the base
has the \emph{weak} Hilbert property and the fibres have the \emph{usual}
Hilbert property:

\begin{theorem}[Javanpeykar]\label{Thm:JavMixedFibration}
    Let $k$ be a number field
    and let $X \to S$ be a morphism of normal integral projective $k$-varieties.
    Let $\Omega\subseteq S(k)$ be a subset that is not strongly thin in $S$ and such that,
    for every $s \in \Omega$, the fibre $X_s$ is a normal integral variety satisfying the Hilbert property over $k$.
    Then $X$ satisfies the weak Hilbert property over $k$.
\end{theorem}

The purpose of this short note is to extend the above ``mixed fibration theorem'' to not necessarily proper varieties over an arbitrary field of characteristic zero:

\begin{theorem}\label{Thm:MixedFibration}
    Let $k$ be a field of characteristic zero and let $X \to S$ be a dominant morphism of normal integral varieties over $k$.
    Let $\Gamma \subseteq X(k)$ and $\Sigma \subseteq S(k)$ be subsets such that
    $\Sigma$ is not strongly thin in $S$ and, for every $s \in \Sigma$, the fibre $X_s$ is a normal integral $k$-variety and $\Gamma \cap X_s(k)$ is not thin in $X_s$.
    Then $\Gamma$ is not strongly thin in $X$.
\end{theorem}

The proof of Theorem~\ref{Thm:JavMixedFibration} uses Corvaja--Zannier's result \cite[Theorem~1.4]{CZHP} that a normal integral projective variety with the Hilbert property over a number field is algebraically simply connected
in order to conclude that the covers of the fibres $X_s$ occurring in the proof are ramified, and thus of degree at least $2$.
We avoid this argument by comparing the degrees of these covers of $X_s$ to the degrees of other non-trivial covers,
and by using Nagata compactification to get rid of the projectivity assumption.

Our motivation for writing this short note is that the mixed fibration in the generality proven here allows for new applications to the conjectures of Campana and Corvaja--Zannier
(see \cite[Lemma~7.7]{BartschJavanpeykarLevin} or \cite[Corollary~7.8]{BartschJavanpeykarLevin} for concrete examples).

\begin{proof}[Proof of Theorem~\ref{Thm:MixedFibration}]
    Let $(\pi_i \colon Y_i \to X)_{i=1}^n$
    be a finite collection of ramified covers.
    To prove the claim, we have to show that $\Gamma \setminus \cup_{i=1}^n \pi_i(Y_i(k))$ is dense in $X$.
    Let $X \to \bar X \to S$ be a factorization of $X \to S$ such that $\bar X$ is a normal integral $k$-variety, $X \to \bar X$ is an open immersion, and $\bar X \to S$ is proper.
    (This can be achieved by taking a Nagata compactification of $X \to S$ \cite[Tag~0F3T]{stacks-project} and normalization.)
    For every $i=1,\ldots,n$, let $\bar \pi_i \colon \bar Y_i \to \bar X$ be the normalization of $\bar X$ in $Y_i$
    and let $\bar Y_i \to T_i \overset{\psi_i}{\longrightarrow} S$ be the Stein factorization of $\bar Y_i \to S$ (see \cite[Tag~03H0]{stacks-project}), so $\psi_i$ is finite, the fibres of $\bar Y_i \to T_i$ are geometrically connected, and $T_i$ is the normalization of $S$ in $\bar Y_i$.
    We have the following commutative diagram.
    \[
        \xymatrix{
            Y_i \ar_{\pi_i}[d] \ar[r] & \bar Y_i \ar_{\bar \pi_i}[d] \ar[dr] \\
            X \ar[dr] \ar[r] & \bar X \ar[d] & T_i \ar^{\psi_i}[dl] \\
            & S
        }
    \]
    Since $\bar Y_i$ is normal and integral, it follows from \cite[Tag~035L]{stacks-project} that $T_i$ is normal and integral.

    Let $U \subseteq X$ be a dense open subscheme such that,
    for every $i$, the restriction
    $Y_{i,U} := Y_i \times_X U \to U$ is étale.
    Moreover, let
    $V \subseteq S$ be a dense open such that, for every $s \in V$, we have $U_s \subseteq X_s$ dense open.
    (Such $V$ exists by \cite[Tag~0573]{stacks-project} since $X$ and the fibres $X_s$ are reduced, which implies that scheme-theoretic density and Zariski-density of open subschemes are equivalent by \cite[Tag~056D]{stacks-project}.)

    Up to reindexing, we may and do assume that there exists an integer $r \in \{1,\ldots,n\}$
    such that $\psi_i$ is unramified (hence étale) for $i=1,\ldots,r$ and ramified for $i=r+1,\ldots,n$.
    Let $V' \subseteq S$ be a dense open subscheme such that,
    for every $i \leq r$ and every $s\in V'$, the fibre $\bar Y_{i,s}$ is normal.
    Define
    \[
        \Omega = (\Sigma \cap V \cap V') \setminus \bigcup_{i=r+1}^n \psi_i (T_i(k))
    .\]
    Since $\Sigma$ is not strongly thin in $S$ and $V, V'$ are dense opens, $\Omega$ is dense in $S$.
    Moreover, for every $s \in \Omega$ and every $i \geq r+1$,
    we have $Y_{i,s}(k) = \emptyset$ (as $T_{i,s}(k)$ is empty).
    Let
    \[
        \Psi = \bigcup_{s \in \Omega} \left( (\Gamma \cap X_s)
        \setminus \cup_{i=1}^r \pi_{i,s} (Y_{i,s}(k)) \right)
        \subseteq \Gamma \setminus \bigcup_{i=1}^n \pi_i(Y_i(k))
    .\]
    Since $\Omega$ is dense in $S$, to conclude the proof, it suffices to show that, for every $s \in \Omega$, the set
    $(\Gamma \cap X_s) \setminus \cup_{i=1}^r \pi_{i,s} (Y_{i,s}(k))$ is dense in $X_s$.

    Let $i\leq r$ and $s\in \Omega$. By the étaleness of $\psi_i$,
    the fibre $T_{i,s}$ is a disjoint union of closed points,
    say $T_{i,s} = t_1 \sqcup \ldots \sqcup t_{N}$.
    We then have $\bar Y_{i,s} = \sqcup_{j=1}^{N} \bar Y_{i,t_j}$ with $\bar Y_{i,t_j}$ normal and geometrically connected, and thus geometrically integral, over the residue field $\kappa(t_j)$.
    Moreover, we have $Y_{i,s} = \sqcup_{j=1}^{N} Y_{i,t_j}$ and, for every $j$,
    the fibre $Y_{i,t_j}$ is either empty or a dense open of $\bar Y_{i,t_j}$ (in which case it is normal and geometrically integral over $\kappa(t_j)$).
    Clearly, if $\kappa(t_j) \neq k$, then $Y_{i,t_j}(k) = \emptyset$.
    Therefore, since $\Gamma \cap X_s$ is not thin in $X_s$, it suffices to show that, for every~$j$
    with $\kappa(t_j) = k$ and $Y_{i,t_j} \subseteq \bar Y_{i,t_j}$ dense open,
    the degree of $Y_{i,t_j} \to X_s$ is at least~$2$.

    Let $X_{T_i} = X \times_S T_i$ and $Y_{i,T_i} = Y_i \times_S T_i$.
    We have $T_i \to S$ finite étale, so $X_{T_i} \to X$ is finite étale.
    It follows that the induced morphism $Y_i \to X_{T_i}$ is not étale (since $Y_i \to X$ would be étale otherwise).
    Since $X_{T_i} \to X$ and $Y_i \to X$ are finite, $Y_i \to X_{T_i}$ is also finite \cite[Tag~035D]{stacks-project}.
    Moreover, $Y_i \to X_{T_i}$ is a left factor of the surjective morphism $Y_{i,T_i} \to X_{T_i}$ and therefore surjective, so $X_{T_i}$ is connected.
    By Zariski's main theorem \cite[Tag~0AB1]{stacks-project} and the fact that $Y_i \to X_{T_i}$ is not étale (in particular not an isomorphism), it follows that $Y_i \to X_{T_i}$ is of degree at least $2$.
    
    Let $t_j \in T_{i,s}(k)$ such that $Y_{i,t_j}\neq \emptyset$.
    To conclude the proof, we show that the degree of $Y_{i,t_j} \to X_s$ is at least $2$.
    Note that the composition of $t_j \colon \Spec k \to T_i$ and $\psi_i$ is $s$,
    so we have $X_{T_i,t_j} = X_s$.
    Let $U_{T_i} = U \times_S T_i$.
    Then $Y_{i,U} \subseteq Y_i$ and $U_{T_i} \subseteq X_{T_i}$ are dense opens
    and the induced morphism $Y_{i,U} \to U_{T_i}$ has the same degree as $Y_i \to X_{T_i}$, which is at least $2$.
    By the definition of $U$, the morphism $Y_{i,U} \to U$ is étale, and since $U_{T_i} \to U$ is étale, it follows from \cite[Tag~02GW]{stacks-project} that $Y_{i,U} \to U_{T_i}$ is étale, hence flat.
    By flatness, the degree of the induced morphism on fibres $(Y_{i,U})_{t_j} \to U_s$ 
    equals the degree of $Y_{i,U} \to U_{T_i}$, i.e., is at least $2$.
    Since $(Y_{i,U})_{t_j} \to U_s$ is the pullback of $Y_{i,t_j} \to X_s$ along $U_s \to X_S$ and $U_s \subseteq X_s$ is a dense open,
    this shows that $Y_{i,t_j} \to X_s$ is of degree at least $2$, as required.
\end{proof}

\begin{ack}
We thank Ariyan Javanpeykar for suggesting to us to write this note and for helpful discussions.
We also thank the referee for helpful comments.
This research was funded by the Studienstiftung des Deutschen Volkes (German Academic Scholarship Foundation). We gratefully acknowledge support by the Deutsche Forschungsgemeinschaft (DFG, German Research Foundation) – Project-ID 444845124 – TRR 326.     
\end{ack}

\bibliography{references}{}
\bibliographystyle{alpha}

\end{document}